\documentclass[12pt]{article}
\usepackage[utf8]{inputenc}
\usepackage{graphicx}
\usepackage{amsmath, amssymb,theorem,verbatim}
\usepackage{hyperref}
\usepackage{textcomp} 
\usepackage{setspace}
\usepackage{mathtools}
\usepackage[round]{natbib}
\usepackage{multirow} 
\usepackage{authblk}
\usepackage[usenames,dvipsnames]{xcolor}

\numberwithin{equation}{section}

\newtheorem{theorem}{Theorem}[section]
\newtheorem{lemma}{Lemma}[section] 
\newtheorem{assumption}{Assumption}[section] 
 
\newtheorem{corollary}{Corollary}[section] 

\newtheorem{remark}{Remark}[section] 

\newcommand{\cov}{\mathrm{cov}}
\newcommand{\spa}{\mathrm{sp}}
\newcommand{\var}{\mathrm{var}}

\newcommand{\Ex}{\mathbb{E}}

\title{A prediction perspective on the Wiener-Hopf equations for time series}

\author[1]{Suhasini Subba Rao \thanks{\texttt{suhasini@stat.tamu.edu}}
}
\author[2]{Junho Yang \thanks{\texttt{junhoyang@stat.sinica.edu.tw}, Authors ordered alphabetically.}
}
\affil[1]{Texas A\&M University, College Station, TX 77845,  U.S.A.}
\affil[2]{Academia Sinica, Taipei 115, Taiwan}

\date{\today}

\begin{document}

\maketitle

\begin{abstract}
The Wiener-Hopf equations are a Toeplitz system of linear equations
that naturally arise in several applications in time series. These include the
update and prediction step of the stationary 
Kalman filter equations and the prediction of bivariate time
series. The celebrated Wiener-Hopf technique is usually used for solving
these equations and is based on a comparison of coefficients in a
Fourier series expansion. However, a statistical interpretation of
both the method and solution is opaque. 
The purpose of this note is to revisit the (discrete) Wiener-Hopf
equations and obtain an alternative solution that is more aligned with
classical techniques in time series analysis. 
Specifically, we propose a solution to the Wiener-Hopf equations that
combines linear prediction with deconvolution. 

The Wiener-Hopf solution requires the
spectral factorization of the underlying spectral density
function. For ease of evaluation it is often assumed that the spectral
density is rational. This  allows one to obtain a computationally
tractable solution.
However, this leads to an approximation error when the
underlying spectral density is not a rational function. We use the
proposed solution with Baxter's inequality to derive an error bound for the rational spectral density approximation.

\vspace{0.5em}


\end{abstract}

\section{Introduction} \label{sec:intro}

The Wiener-Hopf technique (\cite{p:hop-wie-31, b:hop-34}) was first proposed in the 1930s as a method for solving an integral equation of the form 
\begin{equation*}
g(\tau) = \int_{0}^{\infty}h(t)c_{}(\tau-t)dt \qquad \text{for} \qquad \tau \in [0, \infty)
\end{equation*}
in terms of $h(\cdot)$, where $c(\cdot)$ is a known difference kernel and $g(\cdot)$ is a specified
function. The above integral equation and the Wiener-Hopf
technique have been widely used in many applications
in applied mathematics and engineering (see \cite{p:law-07} for a review).
In the 1940s, \cite{b:wie-49} reformulated the problem within 
 discrete time, which is commonly referred to as the Wiener (causal) filter.
The discretization elegantly encapsulates several problems in time series analysis. For example,
the best fitting finite order autoregressive parameters
fall under this framework. The autoregressive parameters can be expressed as a solution of a system of finite interval Wiener-Hopf equations
 (commonly referred to as the FIR Wiener filter), for which \cite{p:lev-47} and \cite{p:dur-60} proposed a $O(n^{2})$ method for solving these equations. 
More broadly, the best linear predictor of a causal stationary time series naturally gives rise
to the Wiener filter, e.g., the prediction of hidden states 
in a Kalman filter model. 
The purpose of this paper is to revisit the discrete-time Wiener-Hopf equations (it is precisely defined in (\ref{eq:1sideconvol}))  and derive an alternative solution using the tools of linear prediction. 
Below we briefly review some classical results on the Wiener filter.

Suppose that $\{X_t: t\in \mathbb{Z}\}$ is a real-valued, zero mean weakly stationary time series defined on the probability space $(\Omega, \mathcal{F}, P)$ and $c(r) = \cov(X_0, X_{-r})$ is the autocovariance function of $\{X_t\}$.
Let $\mathcal{H}_{\infty}$ and $\mathcal{H}_{t}$ ($t\in \mathbb{Z}$) denote closed sub-spaces of the real Hilbert space $L_2(\Omega, \mathcal{F}, P)$ 
spanned by $\{X_t: t \in \mathbb{Z}\}$ and $\{X_j: j \leq t\}$ respectively. We denote the
orthogonal projection onto the closed subspace $V \in L_2(\Omega, \mathcal{F},
P)$ as $P_{V}$. For $Y \in L_2(\Omega, \mathcal{F}, P)$, 
the orthogonal projection of $Y$ onto $\mathcal{H}_{0}$ is  
\begin{equation}
\label{eq:PY1}
P_{\mathcal{H}_{0}}(Y) =  \sum_{j=0}^{\infty}h_{j}X_{-j},
\end{equation}
where $P_{\mathcal{H}_{0}}(Y)=\arg\min_{U\in \mathcal{H}_0} \Ex|Y_{0} - U|^2$. 
To evaluate $\{h_{j}: j \geq 0\}$, we rewrite (\ref{eq:PY1}) as a system of
normal equations. By using that $P_{\mathcal{H}_{0}}(Y_{})$ is an
orthogonal projection onto $\mathcal{H}_{0}$
it is easily shown that (\ref{eq:PY1}) leads to
the system of normal equations
\begin{equation} 
\label{eq:1sideconvol}
c_{YX}(\ell) =  \sum_{j=0}^{\infty}h_{j}c_{}(\ell-j) \qquad \text{for} \qquad \ell \geq 0,
\end{equation} 
where $c_{YX}(\ell) = \cov(Y_{},X_{-\ell})$. 
The above set of equations is typically referred to as the discrete-time
Wiener-Hopf equations (or semi-infinite Toeplitz equations).  
There are two well-known methods for solving this equation in the frequency domain; 
the Wiener-Hopf technique (sometimes called the gapped function, see \cite{b:wie-49})
and the prewhitening method proposed by \cite{p:bod-sha-50} and
\cite{p:zad-rag-50}.
Both solutions solve 
for $H(\omega) = \sum_{j=0}^{\infty} h_je^{ij\omega}$ (see
\cite{p:kal-74}, \cite{b:kal-80} and \cite{b:orf-18}, Sections 11.3-11.8).
The Wiener-Hopf technique is based on the spectral factorization and a comparison of Fourier coefficients corresponding to the 
negative and non-negative indices in Fourier series expansion.
The prewhitening method, as the name suggests, is more in the spirit of time series where the time series $\{X_{t}\}$ is ``whitened'' using an autoregressive filter. 

To state the solution, we assume the spectral density $f(\omega) =
\sum_{r\in\mathbb{Z}} c(r)e^{ir\omega}$  satisfies the condition
$0< \inf_{\omega}f_{}(\omega) \leq \sup_{\omega}f(\omega) < \infty$. 
Then, 
$\{X_t\}$ admits an infinite order Wold-type MA and AR representation (\cite{b:pou-01}, Sections 5-6 and 
\cite{p:kra-18}, page 706)
\begin{equation} 
\label{eq:Wold}
X_t = \varepsilon_t + \sum_{j=1}^{\infty} \psi_j \varepsilon_{t-j}, \qquad
X_t - \sum_{j=1}^{\infty} \phi_j X_{t-j} = \varepsilon_t \qquad t\in \mathbb{Z}
\end{equation} 
where $\sum_{j=1}^{\infty} \psi_j^2 <\infty$, $\sum_{j=1}^{\infty} \phi_j^2 <\infty$, and 
$\{\varepsilon_t\}$ is a uniquely determined white noise process 
with $\Ex \varepsilon_t^2
= \sigma^2 >0$ and $\varepsilon_t$ is orthogonal to $\mathcal{H}_{t-1}$ for $t\in \mathbb{Z}$.
Insights into how the Wold representation in (\ref{eq:Wold})
connects linear prediction through the infinite order MA and AR coefficients is given in \cite{p:che-93} and \cite{p:mey-15}. We mention that (\ref{eq:Wold}) holds under the weaker condition that  $\inf_\omega f(\omega) >0$ (see, for example, \cite{p:wie-mas-58}).

From (\ref{eq:Wold}), we immediately obtain the spectral factorization
$f_{}(\omega) = \sigma^{2}|\phi_{}(\omega)|^{-2}$,
where 
$\phi(\omega) =  1 - \sum_{j=1}^{\infty} \phi_j
e^{ij\omega}$. 
Given $A(\omega) =
\sum_{j=-\infty}^{\infty}a_{j}e^{ij\omega}$, we use the notation
 $[A(\omega)]_{+} = \sum_{j=0}^{\infty}a_{j}e^{ij\omega}$ and 
$[A(\omega)]_{-} = \sum_{j=-\infty}^{-1}a_{j}e^{ij\omega}$. 
Both the Wiener-Hopf technique and prewhitening method yield the solution
\begin{equation}
\label{eq:GCov4}
H(\omega)  
= \sigma^{-2} \phi_{}(\omega)[\phi(\omega)^{*}
  f_{YX}(\omega)]_{+},
\end{equation}
where $f_{YX}(\omega)  =
\sum_{\ell\in \mathbb{Z}}c_{YX}(\ell)e^{i\ell \omega}$ and
$\phi(\omega)^*$ is a complex conjugate of $\phi(\omega)$. 

We mention that the special case of $m$-step ahead forecasts falls
under this framework. By setting $Y_{} = X_{m}$, the coefficients
$\{h_j: j\geq0\}$ are the $m$-step ahead prediction coefficients and the solution for 
$H(\omega)$ is
\begin{equation*}
H_{m}(\omega) = \phi(\omega) [\psi(\omega)e^{-im \omega}]_+ \qquad m >0,
\end{equation*} where $\psi(\omega) = (\phi(\omega))^{-1} = 1+\sum_{j=1}^{\infty} \psi_je^{ij\omega}$
and MA coefficients $\{\psi_j\}$ is from (\ref{eq:Wold}).

The normal equations in (\ref{eq:1sideconvol}) belong to the general
class of Wiener-Hopf equations of the form 
\begin{equation}
\label{eq:gell}
g_{\ell} =  \sum_{j=0}^{\infty}h_{j}c(\ell-j) \qquad \text{for} \qquad \ell \geq 0,
\end{equation} 
where $\{c(r): r\in \mathbb{Z}\}$ is a symmetric, positive
definite sequence.  The Wiener-Hopf technique yields the solution
\begin{equation}
\label{eq:GCov42}
H(\omega)  = 
\sigma^{-2}\phi(\omega)[\phi_{}(\omega)^{*} G_+(\omega)]_{+},
\end{equation}
where $G_+(\omega) =
\sum_{\ell=0}^{\infty}g_{\ell}e^{i\ell\omega}$ (the derivation is well-known,
but for completeness we give  a short proof in Section \ref{sec:general}). 
An alternative method for solving for $\{h_{j}: j \geq 0\}$ is
within the time domain.  
This is done by representing (\ref{eq:gell}) as the semi-infinite Toeplitz system
\begin{equation}
\label{eq:gell1}
{\bf g}_+ = T(f){\bf h}_+, 
\end{equation}
where ${\bf g}_+ = (g_{0},g_{1},g_{2},\ldots)^{\prime}$ and ${\bf h}_+ =
(h_{0},h_{1},\ldots)^{\prime}$ are semi-infinite column vectors (sequences)
and $T(f)$ is a Toeplitz matrix of the form 
$T(f) = (c(t-\tau); t,\tau\geq 0)$. Let $\{ \phi_{j}: j \geq 0\}$ (setting $\phi_0 = -1$) denote the infinite order AR coefficients corresponding to $f$ defined as in (\ref{eq:Wold}) and
$\phi(\cdot)$ be its Fourier transform. 
By letting $\phi_{j}=0$ for $j <0$, we define the
lower triangular Toeplitz matrix $T(\phi) = (\phi_{t-\tau};t,\tau\geq
0)$. Provided that $0<\inf_\omega f(\omega) \leq \sup_\omega f(\omega)
<\infty$, it is well-known that $T(f)$ is invertible on 
$\ell_2^+ = \{(v_0, v_1, ...)^{\prime}: \sum_{j=0}^{\infty} |v_j|^2<\infty\}$, 
and the inverse is $T(f)^{-1} = \sigma^{-2}T(\phi)T(\phi)^{*}$ 
(see, for example, Theorem III of \cite{p:wid-60}). Thus, the time domain solution to (\ref{eq:gell}) is 
${\bf h}_+ = T(f)^{-1}{\bf g}_+ = \sigma^{-2}T(\phi)T(\phi)^{*}{\bf g}_+$. 

In this paper, we study the Wiener-Hopf equations from a time series perspective, combining 
the prediction theory developed in the time domain with the deconvolution method
in the frequency domain. Observe that (\ref{eq:gell})
is a system of semi-infinite convolution equations (since the
equations only hold for non-negative index $\ell$), thus the
standard deconvolution approach  is not possible. 
In \cite{p:sub-yang-21}, we used the tools of linear prediction to
rewrite the Gaussian likelihood of a stationary time
series within the frequency domain. We transfer some of these ideas to
solving the Wiener-Hopf equations.  In Section \ref{sec:prediction}, we show that
we can circumvent the constraint $\ell \geq 0$, by using linear prediction to yield
the normal equations in (\ref{eq:1sideconvol}) for all $\ell\in
\mathbb{Z}$. In Section \ref{sec:general}, 
we show that there exists a stationary time series $\{X_{t}\}$ and random variable $Y_{} \in \mathcal{H}_{0}$ where $Y_{}$ and $\{X_t\}$
induce the general Wiener-Hopf equations of the form (\ref{eq:gell}).
This allows us to use the aforementioned technique to reformulate the Wiener-Hopf equations as a
bi-infinite Toeplitz system, and thus obtain a solution to $H(\omega)$ as a deconvolution. 
The same technique is used to obtain an expression for entries of the inverse Toeplitz
matrix $T(f)^{-1}$.

In practice, evaluating $H(\omega)$ in (\ref{eq:GCov4}) requires 
the spectral factorization of the underlying spectral density.
One strategy is to assume that the spectral density is rational, which
allows one to obtain a 
computationally tractable solution for $H(\omega)$. 
Of course, this leads to an approximation error in $H(\omega)$ when the
underlying spectral density is not a rational function. 
In Section \ref{sec:approx}, we show that Baxter's inequality (\cite{p:bax-62, p:bax-63}) can be utilized to obtain a bound between
$H(\omega)$ and its approximation based on a rational approximation of the general spectral density.
 The proof of the results in Sections \ref{sec:alternative} and \ref{sec:approx} can be found in the Appendix. 


\section{A prediction approach}\label{sec:alternative}

\subsection{Notation and Assumptions}\label{sec:notation}

In this section, we collect together the notation introduced in Section \ref{sec:intro}
and some additional notation necessary for the paper. 

Let $L_{2}([0,2\pi))$ be the space of all square-integrable complex
functions on $[0,2\pi)$ and $\ell_{2}$ 
the space of all bi-infinite complex column vectors ${\bf v} =
(\ldots,v_{-1},v_{0},v_{1},\ldots)^{\prime}$ where $\sum_{j\in
  \mathbb{Z}}|v_{j}|^{2}<\infty$. Similarly, we let 
$\ell_2^+ = \{ {\bf v}_+ = (v_0, v_1, ...)^{\prime}:
\sum_{j=0}^{\infty}|v_j|^2 <\infty \}$  denote the space of all
semi-infinite square summable vector sequences. 
To connect the time and frequency domain through an isomorphism, we define the Fourier transform 
$F:\ell_{2} \rightarrow L_{2}([0,2\pi))$ 
\begin{equation*}
F({\bf v})(\omega) = \sum_{j\in \mathbb{Z}}v_{j}e^{ij\omega}.
\end{equation*}
For $f(\omega) = \sum_{r\in \mathbb{Z}}c(r)e^{ir\omega} \in L_{2}([0,2\pi))$, define the semi- and bi-infinite Toeplitz operators $T(f)$ and $T_\pm (f)$ 
on $\ell_2^+$ and $\ell_2$ with the matrix form
$T(f) = (c(t-\tau);t,\tau\geq 0)$ and 
 $T_{\pm}(f) = (c(t-\tau);t,\tau \in \mathbb{Z})$, respectively.
This paper will make frequent use of the convolution
 theorem: If ${\bf h} \in \ell_{2}$, then
 $F(T_{\pm}(f){\bf h})(\omega) = f(\omega) F({\bf h})(\omega)$.


\begin{assumption}
\label{assum:A} Let $\{c(r):r \in \mathbb{Z}\}$ be a
symmetric positive definite sequence on $\ell_2$ and $f(\omega) = \sum_{r\in
  \mathbb{Z}} c(r)e^{ir\omega}$ be its Fourier transform. 
Then,
\begin{itemize}
\item[(i)] $0<\inf_{\omega}f_{}(\omega) \leq \sup_{\omega} f(\omega) <\infty$.

\item[(ii)] For some $K>1$ we have $\sum_{r \in \mathbb{Z}}|r^{K}c(r)|<\infty$.
\end{itemize}
\end{assumption}
Under Assumption \ref{assum:A}(i), we have the unique factorization
\begin{equation}
\label{eq:phipsidef}
f_{}(\omega) = \sigma^{2}|\psi_{}(\omega)|^{2} =
\sigma^{2}|\phi_{}(\omega)|^{-2},
\end{equation} 
where $\sigma^{2}>0$,
$\psi_{}(\omega) = 1+\sum_{j=1}^{\infty}\psi_{j}e^{ij\omega}$ and
 $\phi_{}(\omega) = (\psi_{}(\omega))^{-1} =
 1-\sum_{j=1}^{\infty}\phi_{j}e^{ij\omega}$. 
The characteristic polynomials  $\Psi(z) = 1 + \sum_{j=1}^{\infty} \psi_j z^j$ and $\Phi(z) =
1-\sum_{j=1}^{\infty} \phi_j z^j$ do not have zeroes in 
$|z|\leq 1$ thus the AR$(\infty)$ parameters are causal or
equivalently are said to have minimum phase
 (see \cite{p:sze-21} and \cite{p:ino-00}, pages 68-69).

We mention that Assumption \ref{assum:A}(i) is used in all the results
in this paper, whereas Assumption \ref{assum:A}(ii) is only required
in the approximation theorem in Section \ref{sec:approx}. Under Assumption \ref{assum:A}(ii), 
both $\sum_{j=1}^{\infty} j^{K} |\psi_j|$ and $\sum_{j=1}^{\infty} j^{K}|\phi_j|$ are finite (see \cite{p:che-93} and \cite{p:mey-15}).

\subsection{Solving Wiener-Hopf equations using linear prediction}\label{sec:prediction}

We now give an alternative formulation for the solution of
(\ref{eq:1sideconvol}) and (\ref{eq:gell}), which 
utilizes properties of linear prediction to solve it using a standard deconvolution method. 
To integrate our derivation within the Wiener causal filter framework, we start with the classical Wiener filter.
For $Y \in L_2(\Omega, \mathcal{F}, P)$, consider the projection of $Y$ onto $\mathcal{H}_0$
\begin{equation}
\label{eq:PY}
P_{\mathcal{H}_{0}}(Y_{}) =  \sum_{j=0}^{\infty}h_{j}X_{-j}.
\end{equation} 
We observe that by construction, (\ref{eq:PY}) gives rise to the normal equations 
\begin{equation}
\label{eq:1sideconvol1}
\cov(Y,X_{-\ell})   = \sum_{j=0}^{\infty}h_{j}c(\ell-j) \qquad \ell \geq 0.
\end{equation}
Since (\ref{eq:1sideconvol1}) only holds for non-negative
$\ell$, this prevents one using 
 deconvolution to solve for $H(\omega)$. Instead, we define a
``proxy'' set of variables for $\{X_{-\ell}: \ell <0 \}$ such that
(\ref{eq:1sideconvol1}) is valid for $\ell<0$. 
By using the property of orthogonal projections, we have
\begin{equation*}
\cov(Y, P_{\mathcal{H}_0}(X_{-\ell}))=
 \cov( P_{\mathcal{H}_{0}}(Y_{}), X_{-\ell}) \qquad \ell < 0.
\end{equation*}
This gives 
\begin{equation} 
\label{eq:1sideconvol2}
\cov(Y,P_{\mathcal{H}_0}(X_{-\ell}) )
=  \sum_{j=0}^{\infty}h_{j}\cov (X_{-j},X_{-\ell})
=  \sum_{j=0}^{\infty}h_{j} c_{}(\ell-j) \qquad \ell < 0.
\end{equation}
Equations (\ref{eq:1sideconvol1}) and (\ref{eq:1sideconvol2})  allow
us to represent the solution of $H(\omega)$ as a deconvolution. 
We define the semi- and bi-infinite sequences
 ${\bf c}_{-} = (\cov(Y,P_{\mathcal{H}_0}(X_{-\ell}));\ell <0)^{\prime}$,
${\bf c}_{+}=(\cov(Y,X_{-\ell});\ell \geq 0)^{\prime}$, and
  ${\bf c}_{\pm} = ({\bf c}_{-}^{\prime}, {\bf c}_{+}^{\prime})^\prime$. Taking the Fourier
  transform of ${\bf c}_{\pm}$ and using the convolution theorem gives
  $F({\bf c}_{\pm})(\omega) = H(\omega)f(\omega)$. Thus 
\begin{equation}
\label{eq:HFF}
H(\omega) = \frac{F({\bf c}_{\pm})(\omega)}{f(\omega)}
=
\frac{\sum_{\ell = 0}^{\infty} \cov (Y, X_{-\ell})e^{i\ell\omega} 
+\sum_{\ell=1}^{\infty}\cov\left( Y,P_{\mathcal{H}_0}(X_{\ell})
              \right)e^{-i\ell\omega}}{f(\omega)}.
\end{equation} 
This forms the key to the following theorem. 

\begin{theorem}\label{theorem:wiener}
Suppose $\{X_{t}\}$ is a stationary time series on the probability space $(\Omega, \mathcal{F}, P)$
whose spectral density satisfies Assumption \ref{assum:A}(i). 
Let $\phi_{}(\cdot)$ and
$\psi_{}(\cdot)$ are be defined as in (\ref{eq:phipsidef}) and 
$\phi_{\ell}(\omega) = \sum_{s=1}^{\infty}\phi_{\ell+s}e^{is\omega}$ for $\ell \geq 0$.
For $Y \in L_2(\Omega, \mathcal{F}, P)$, let $P_{\mathcal{H}_{0}}(Y_{}) =  \sum_{j=0}^{\infty}h_{j}X_{-j}$. Then, $(h_{j};j \geq 0)^{\prime} \in \ell_{2}^{+}$, $( c_{YX}(\ell) = \cov (Y, X_{-\ell}); \ell \geq 0)^{\prime} \in
\ell_2^+$ and 
\begin{equation}
\label{eq:Homega}
H(\omega) 
=\frac{\sum_{\ell =  0}^{\infty}c_{YX}(\ell) \big( e^{i\ell\omega} +
    \psi_{}(\omega)^{*}\phi_{\ell}(\omega)^{*}\big)}{f_{}(\omega)}.
\end{equation} 
The above solution can alternatively be expressed as 
\begin{equation}
\label{eq:alternativeform}
H(\omega)=\sigma^{-2} \phi(\omega)
\sum_{\ell =0}^{\infty}c_{YX}(\ell) \bigg( e^{i\ell\omega}-\sum_{s=1}^{\ell} \phi_s e^{i(\ell-s)\omega}\bigg),
\end{equation} where for $\ell=0$, $\sum_{s=1}^{0} = 0$.
\end{theorem}
\noindent PROOF. See Appendix \ref{sec:proof}. \hfill $\Box$

\vspace{1em}

We thank an anonymous referee for pointing out that the representation
in (\ref{eq:Homega})  is equivalent to
(\ref{eq:alternativeform}). The
benefit of the latter representation is that
both $\phi(\omega)$ and $e^{i\ell\omega}-\sum_{s=1}^{\ell}
\phi_s e^{i(\ell-s)\omega}$  (for $\ell \geq 0$)  are in terms of the power series of $e^{i\omega}$, thus it is transparent that $H(\omega)$ is causal.

\begin{remark}[Relationship to concurrent filters]
There is a close relationship between Theorem
\ref{theorem:wiener} and solutions to concurrent filters (that are
frequently used by the U.S. Census Bureau). Notable applications are the multi-step ahead forecasts used in the derivation of the X-11 and X-11-ARIMA seasonal filters
(see \cite{p:dag-75, p:dag-82} and \cite{p:lad-que-12} for a review). 
In relation to the Wiener filter, 
this is the technique of using multi-step ahead forecasts to obtain a
solution to concurrent filter $P_{\mathcal{H}_0}(Y) =
\sum_{j=0}^{\infty}h_{j}X_{-j}$ from the two-sided filter 
$P_{\mathcal{H}_\infty}(Y) =\sum_{j=-\infty}^{\infty}a_{j}X_{-j}$
where $\sum_{j=-\infty}^{\infty}a_{j}e^{ij\omega}=\sum_{r\in
  \mathbb{Z}}c_{YX}(r)e^{ir\omega}/f(\omega)$; See
\cite{p:bel-mar-04} and \cite{p:wil-mce-16}, Section 2 (Proposition
1). We summarize the technique below.
By standard projection arguments, we have
\begin{equation} 
\label{eq:PHY}
P_{\mathcal{H}_0}(Y) =
P_{\mathcal{H}_0}P_{\mathcal{H}_\infty}(Y) =
\sum_{j=0}^{\infty} a_{-j} X_{-j} + \sum_{\ell=1}^{\infty} a_{\ell}
P_{\mathcal{H}_0} (X_{\ell}) \
= \sum_{j=0}^{\infty} \left( a_{-j} + \sum_{\ell=1}^{\infty}a_{\ell}\phi_{j}(\ell)
\right) X_{-j},
\end{equation}
where $P_{\mathcal{H}_0} (X_{\ell})=\sum_{j=0}^{\infty} \phi_{j}(\ell)X_{-j}$ (an expression for these
coefficients in terms of AR and MA coefficients is given in Appendix, (\ref{eq:phiellj})).
Therefore, by comparing the above to (\ref{eq:PY}) we have 
 $h_j = a_{-j} + \sum_{\ell=1}^{\infty}a_{\ell} \phi_{j}(\ell)$ for $j
 \geq 0$. Note that both (\ref{eq:alternativeform}) and (\ref{eq:PHY}) 
 yield different solutions to the same Wiener-Hopf equations.
\end{remark}

\begin{remark}[Relationship to prediction]\label{remark:wiener}
It is clear that 
$\sum_{\ell=1}^{\infty}X_{\ell}e^{i\ell\omega}$ is not a well-defined random variable. 
However, it is interesting to note that under Assumption
\ref{assum:A}(ii) (for $K=1$)
$\sum_{\ell=1}^{\infty}P_{\mathcal{H}_0}(X_{\ell})e^{i\ell\omega}$ is
a well defined random variable in $\mathcal{H}_0$ and 
\begin{equation}
\label{eq:PPY}
\sum_{\ell=1}^{\infty}P_{\mathcal{H}_0}(X_{\ell})e^{i\ell\omega} = \psi(\omega)\sum_{j=0}^{\infty} \phi_{j}(\omega) X_{-j}.
\end{equation}
In other words, despite
$\sum_{\ell=1}^{\infty}X_{\ell}e^{i\ell\omega}$ not being well-defined,
informally, its projection onto $\mathcal{H}_0$ does exist. 
\end{remark}

\subsection{General Wiener-Hopf equations}\label{sec:general}

We now generalize the prediction approach in the previous section to general Wiener-Hopf
linear equations which satisfy
\begin{equation}
\label{eq:gplus}
g_{\ell} =  \sum_{j=0}^{\infty}h_{j}c(\ell-j) \qquad \ell \geq 0,
\end{equation}
where $\{g_{\ell}: \ell \geq 0\}$ and $\{c(r): r\in \mathbb{Z}\}$
(which is assumed to be a symmetric, positive definite sequence) are
known. We will obtain a solution similar to (\ref{eq:Homega}) but for the normal equations in (\ref{eq:gplus}). 
We first describe the classical Wiener-Hopf method to solve (\ref{eq:gplus}). 
Since $\{c(r)\}$ is known for all $r\in \mathbb{Z}$,
we extend
(\ref{eq:gplus}) to the negative index $\ell<0$, and define $\{g_{\ell}^{}:\ell <0\}$ as
\begin{equation}
\label{eq:gminus}
g_{\ell}^{} =  \sum_{j=0}^{\infty}h_{j}c(\ell-j) \qquad \textrm{for} \qquad \ell < 0.
\end{equation}
Note that $\{g_{\ell}:\ell <0\}$ is not given, but it 
is completely determined by $\{g_{\ell}: \ell \geq 0\}$ and
$\{c(r)\}$ (this can be seen from (\ref{eq:Gminus}), below).
The Wiener-Hopf technique
evaluates the Fourier transform of the above and isolates the 
non-negative indices in the Fourier series expansion
to yield the solution for $H(\omega)$. Specifically,
evaluating the Fourier transform of (\ref{eq:gplus}) and
(\ref{eq:gminus}) gives 
\begin{equation}
\label{eq:Hf}
f(\omega)H(\omega) =  G_-(\omega) + G_+(\omega)
\end{equation}
where $G_{-}(\omega)=\sum_{\ell = -\infty}^{-1} g_{\ell}e^{i\ell\omega}$ and 
$G_{+}(\omega)=\sum_{\ell = 0}^{\infty} g_{\ell}e^{i\ell\omega}$. 
Replacing
$f(\omega)$ with $\sigma^{2}|\psi(\omega)|^2$ and dividing
the above with $\sigma^{2}\psi(\omega)^{*}$ yields
\begin{equation} 
\label{eq:Hf2}
H(\omega)\psi(\omega) = 
\frac{G_-(\omega)}{\sigma^{2}\psi(\omega)^{*}} + \frac{G_+(\omega)}{\sigma^{2}\psi(\omega)^{*}}
= \sigma^{-2} \phi(\omega)^* G_-(\omega) + \sigma^{-2} \phi(\omega)^* G_+(\omega).
\end{equation}
Isolating the  non-negative indices in (\ref{eq:Hf2}) gives the solution
\begin{equation}
\label{eq:HHPP}
H(\omega) = \sigma^{-2}\phi(\omega) 
[ \phi(\omega)^{*} G_+(\omega)]_{+},
\end{equation}
this proves the result stated in (\ref{eq:GCov42}).
Similarly, by isolating the negative indices,
we obtain the expression $G_{-}(\omega)=\sum_{\ell = -\infty}^{-1} g_{\ell}e^{i\ell\omega}$ in terms of $f(\omega)$ and $G_{+}(\omega)$
\begin{equation}
\label{eq:Gminus}
G_{-}(\omega) =  -\psi(\omega)^{*} [ \phi(\omega)^{*} G_{+}(\omega)]_{-}.
\end{equation}
Thus (\ref{eq:HHPP}) and (\ref{eq:Gminus}) provide explicit solutions to $H(\omega)$
and $G_{-}(\omega)$ respectively.
However, from a time series perspective, it is difficult to
interpret these formulas. We now obtain an alternative expression for
these solutions based on the linear prediction of random
variables.

We consider the matrix representation, $T(f){\bf h}_{+} ={\bf
  g}_{+}$, in (\ref{eq:gell1}).
We solve $T(f){\bf h}_{+} ={\bf
  g}_{+}$ by embedding
the semi-infinite Toeplitz matrix $T(f)$ on $\ell_2^+$
into the bi-infinite Toeplitz system on $\ell_2$. 
To relate $T(f)$ and $T_{\pm}(f)$ we partition the
bi-infinite Toeplitz matrix $T_{\pm}(f)$ into four sub-matrices
$C_{00} = (c(t-\tau); t,\tau <
0)$, $C_{01} =
(c(t-\tau); t <0, \tau \geq 0)$, $C_{10} =
(c(t-\tau); t \geq 0, \tau < 0)$, and $C_{11} = (c(t-\tau);t,\tau
\geq 0)$. We observe that $C_{11} = T(f)$.
Further, we let ${\bf h}_{\pm} = ({\bf 0}^\prime, {\bf h}_+^\prime)^\prime = (\ldots,0,0,h_{0},h_{1},h_{2},\ldots)^{\prime}$ and ${\bf g}_{\pm} = ({\bf g}_{-}^\prime, {\bf g}_{+}^\prime)^\prime =  (\ldots,g_{-2},g_{-1},g_{0},g_{1},g_{2},\ldots)^{\prime}$ where ${\bf g}_{-} = C_{00}C_{11}^{-1} {\bf g}_{+}$. Then, we obtain the following
bi-infinite Toeplitz system on $\ell_2$
\begin{equation}
\label{eq:TTT}
 T_{\pm }(f){\bf h}_{\pm}=
\begin{pmatrix}
C_{00}  & C_{01} \\
C_{10}   & C_{11} \\
\end{pmatrix}
\begin{pmatrix}
{\bf 0} \\
{\bf h}_{+}
\end{pmatrix}
=
\begin{pmatrix}
 C_{01}{\bf h}_{+}\\
C_{11} {\bf h}_{+}
\end{pmatrix}
=
\begin{pmatrix}
 C_{01}C^{-1}_{11}{\bf g}_{+}\\
{\bf g}_{+} \\
\end{pmatrix} = 
\begin{pmatrix}
{\bf g}_{-}\\
{\bf g}_{+} \\
\end{pmatrix} 
= {\bf g}_{\pm}.
\end{equation}
We note that the non-negative indices in the sequence ${\bf g}_{\pm}$ are $\{g_{\ell}:\ell\geq0\}$, 
but for the negative indices, where $\ell<0$, it is $g_{\ell} = [C_{01}C^{-1}_{11}{\bf  g}_{+}]_{\ell}$
which is identical to $g_{\ell}^{}$ defined in (\ref{eq:gminus}). 
The Fourier transform on both sides in (\ref{eq:TTT}) gives
$f(\omega)H(\omega) = F({\bf g}_{\pm})(\omega)$,
which is identical to (\ref{eq:Hf}).
We now reformulate the above
equation through the lens of prediction. 
To do this, we construct a stationary process $\{X_{t}\}$ and a random
  variable $Y$ on the same probability space which yields 
(\ref{eq:gplus}) as their normal equations. 

We first note that since $\{c(r): r\in \mathbb{Z}\}$ is a symmetric, positive
definite sequence,  there exists a stationary time series $\{X_{t}\}$
with $\{c(r): r\in \mathbb{Z}\}$ as its autocovariance function (see
\cite{b:bro-dav-06}, Theorem 1.5.1). Using this, define the random variable 
\begin{equation}
\label{eq:Ydef}
Y= \sum_{j=0}^{\infty}h_{j}X_{-j}.
\end{equation} 
Provided that ${\bf h}_+ \in \ell_2^+$, then $\Ex[Y^2]<\infty$ and
thus $Y \in \mathcal{H}_{0}$  (we show in Theorem \ref{corollary:wiener}
that this is true if ${\bf g}_+ \in \ell_2^+$).
By (\ref{eq:gplus}), we observe that $\cov(Y,X_{-\ell}) =
\sum_{j=0}^{\infty}h_{j}c(\ell-j)=g_{\ell}$ for all $\ell \geq 0$. 
We now show that for $\ell<0$,
\begin{eqnarray*}
\cov(Y,X_{-\ell}) = [C_{01}C^{-1}_{11}{\bf  g}_{+}]_{\ell} = g_\ell.
\end{eqnarray*}
First, since $Y \in \mathcal{H}_{0}$, then
$\cov(Y,X_{-\ell}) = \cov(P_{\mathcal{H}_0}(Y),X_{-\ell}) = \cov(Y,P_{\mathcal{H}_0}(X_{-\ell}))$.	
Further, for $\ell<0$, the $\ell$th row (where we start the
enumeration of the rows from the bottom) of
$C_{01}C^{-1}_{11}$ contains the coefficients of the best linear predictor of $X_{-\ell}$ given $\mathcal{H}_0$
\begin{equation} 
\label{eq:C01C11}
P_{\mathcal{H}_0}(X_{-\ell}) = \sum_{j=0}^{\infty}[C_{01}C^{-1}_{11}]_{\ell,j}X_{-j} \qquad \ell<0.
\end{equation} 
A detailed calculation of (\ref{eq:C01C11}) is given in the Appendix.
Using the above, we evaluate $\cov(Y,P_{\mathcal{H}_0}(X_{-\ell}))$ for
$\ell<0$
\begin{eqnarray*}
\cov(Y,P_{\mathcal{H}_0}(X_{-\ell})) &=&
  \cov\left(Y,\sum_{j=0}^{\infty}[C_{01}C^{-1}_{11}]_{\ell,j}X_{-j}\right) \\
&=& \sum_{j=0}^{\infty}[C_{01}C^{-1}_{11}]_{\ell,j}\cov(Y,X_{-j})
    \qquad \qquad (\textrm{from (\ref{eq:Ydef}), } g_{j}=\cov(Y,X_{-j})) \\
&=& \sum_{j=0}^{\infty}[C_{01}C^{-1}_{11}]_{\ell,j}g_{j} =[C_{01}C^{-1}_{11}{\bf  g}_{+}]_{\ell}  = g_{\ell}.
\end{eqnarray*}
Thus the entries of ${\bf g}_{\pm}=({\bf g}_{-}^\prime, {\bf g}_{+}^\prime)^\prime $ 
are indeed the covariances:
${\bf g}_{-} 	= (\cov(Y,P_{\mathcal{H}_0}(X_{-\ell}));\ell <0)^\prime$ and 
${\bf g}_{+} =(\cov(Y,X_{-\ell});\ell \geq 0)^\prime$.
This allows us to use Theorem \ref{theorem:wiener} to solve 
general Wiener-Hopf equations. Further, it gives an intuition to (\ref{eq:gminus}) and (\ref{eq:TTT}). 

\begin{theorem}\label{corollary:wiener}
Suppose that $\{c(r): r\in \mathbb{Z}\}$ is a symmetric, positive definite
sequence and its Fourier transform
$f(\omega) = \sum_{r\in \mathbb{Z}}c(r)e^{ir\omega}$ satisfies Assumption
\ref{assum:A}(i). 
We define the (semi) infinite system of equations 
\begin{equation*}
g_{\ell} =
  \sum_{j=0}^{\infty}h_{j}c(\ell-j) \qquad \ell \geq 0,
\end{equation*}
where $(g_{\ell}; \ell \geq 0)^{\prime} \in \ell_2^+$. Then, $(h_j; j\leq 0)^{\prime} \in \ell_2^+$ and 
\begin{equation}
\label{eq:HH}
H(\omega) = 
\frac{\sum_{\ell=0}^{\infty}
g_{\ell} \big( e^{i\ell\omega} +
    \psi_{}(\omega)^*\phi_{\ell}(\omega)^*\big)}{f_{}(\omega)}.
\end{equation}
Moreover, as in Theorem \ref{theorem:wiener}, $H(\omega)$ can be rewritten as
\begin{equation}
\label{eq:HH2}
H(\omega) = \sigma^{-2} \phi(\omega) \sum_{\ell=0}^{\infty}
g_\ell \bigg( e^{i\ell \omega} - \sum_{s=1}^{\ell} \phi_s e^{i(\ell-s)\omega}\bigg).
\end{equation}

\end{theorem}
\noindent PROOF. See Appendix \ref{sec:proof}. \hfill $\Box$

\vspace{1em}

It is interesting to observe that 
the solution for $H(\omega)$ given in  (\ref{eq:HHPP}) was obtained
by comparing the Fourier coefficients, whereas the
solution in Theorem \ref{corollary:wiener} was obtained using linear
prediction. The two solutions are algebraically different.
We now show that they are the same by direct verification. 

\begin{lemma} \label{lemma:verification}
Suppose the same set of assumptions and notation as in Theorem \ref{corollary:wiener}
hold. Then 
\begin{equation}
\label{eq:link42}
[\phi_{}(\omega)^{*}
  G_{+}(\omega)]_{+}= \sum_{\ell=0}^{\infty}
g_\ell \bigg( e^{i\ell \omega} - \sum_{s=1}^{\ell} \phi_s e^{i(\ell-s)\omega}\bigg),
\end{equation}
where $G_{+}(\omega) = \sum_{\ell=0}^{\infty}g_{\ell}e^{i\ell\omega}$.
\end{lemma}


Theorem \ref{corollary:wiener} can be used to obtain an expression for $T(f)^{-1}$. 
As mentioned in Section \ref{sec:intro}, the time domain solution for the inverse Toeplitz matrix is
$T(f)^{-1} = \sigma^{-2}T(\phi)T(\phi)^{*}$. We show below that 
an alternative expression for the entries of $T(f)^{-1} = (d_{k,j};k,j\geq 0)$  can be
deduced using the deconvolution method described in Theorem \ref{corollary:wiener}.

\begin{corollary}\label{theorem:inverse}
Suppose the same set of assumptions and notation as in Theorem \ref{corollary:wiener} hold. Let
${\bf d}_k = (d_{k,j}; j\geq 0)$ denote the $k$th row of $T(f)^{-1}$. 
Then, ${\bf d}_k^\prime \in \ell_2^+$ for all $k\geq 0$ and the Fourier transform
 $D_{k}(\omega) = \sum_{j=0}^{\infty}d_{k,j}e^{ij\omega}$ is 
\begin{equation*}
D_{k}(\omega) 
=\frac{e^{ik\omega} +
    \psi_{}(\omega)^*\phi_{k}(\omega)^*}{f(\omega)} =
\sigma^{-2} \phi(\omega)\bigg( e^{ik\omega} - \sum_{s=1}^{k} \phi_s e^{i(k-s)\omega}\bigg) \qquad k\geq 0.
\end{equation*}
Therefore,
\begin{eqnarray} 
d_{k,j}  
&=& \frac{\sigma^{-2}}{2\pi}\int_{0}^{2\pi}\phi(\omega)\bigg(
    e^{ik\omega} - 
\sum_{s=1}^{k} \phi_s e^{i(k-s)\omega}\bigg) e^{-ij\omega}d\omega
\qquad j,k\geq 0.
\label{eq:sij} 
\end{eqnarray}
\end{corollary}
\noindent PROOF. See Appendix \ref{sec:proof}. \hfill $\Box$

\begin{remark}[Connection to the inverse of finite order Toeplitz
  matrix]
Consider the $n\times n$ Toeplitz matrix $T_n(f) = (c(s-t);
0 \leq s,t \leq n-1)$ and $d_{k,j}^{(n)} =(T_n(f)^{-1})_{k,j}$.  There
are several different expressions for $d_{k,j}^{(n)}$ including the
Cholesky decomposition given in \cite{p:aka-69, b:pou-01}, and \cite{p:jen-21} or expressions based on a dual process
representation; \cite{p:sub-yang-21} and \cite{p:ino-21}. 
The arguments in this paper can also be used to obtain an alternative expression for the
inverse of a finite dimensional Toeplitz matrix. 
Using similar arguments to those used to prove Corollary
\ref{theorem:inverse}, we obtain
\begin{equation}
\label{eq:dkjn}
d_{k,j}^{(n)} = \frac{1}{2\pi}\int_{0}^{2\pi}\sum_{\ell\in \mathbb{Z}}\frac{\phi_{k}^{(n)}(\ell)
e^{i\ell\omega}}{f(\omega)}e^{-ij\omega}d\omega = \sum_{\ell\in \mathbb{Z}}\phi_{k}^{(n)}(\ell)\gamma(j-\ell) \qquad 0\leq j,k \leq n-1,
\end{equation}
where  $\gamma(k)=(2\pi)^{-1}\int_{0}^{2\pi}f(\omega)^{-1}e^{-ik\omega}d\omega$
(usually called the inverse autocovariance function) and
$\{\phi_{s}^{(n)}(\ell)\}_{s=0}^{n-1}$ are the multi-step ahead finite prediction coefficients;
$P_{\mathcal{H}_{[-(n-1),0]}}(X_{\ell}) = \sum_{s=0}^{n-1}\phi_{s}^{(n)}(\ell)X_{-s}$, where $\mathcal{H}_{[-(n-1),0]} = \spa \{X_{t}: -(n-1) \leq t \leq 0\}$.

It is interesting to compare and contrast (\ref{eq:dkjn}) with the
entries of the finite dimension Cholesky decomposition $T_{n}(f)^{-1} = L_{n}(\phi)L_{n}(\phi)^{\top}$
(see the aforementioned references). Equation 
(\ref{eq:dkjn}) is in terms of products of
coefficients of finite predictors for $X_{\ell}$ ``outside'' the
interval $\{-(n-1),\ldots,0\}$, while 
 the Cholesky decomposition is based on the coefficients of the best linear predictor
of $X_{\ell}$  ``inside'' the interval $\{-(n-1),\ldots,0\}$. 
\end{remark}

\begin{remark}[Multivariate extension]
 The case that the (autocovariance)
sequence $\{{\bf C}(r):r\in \mathbb{Z}\}$ is made up of $d\times d$-dimensions, has not been considered in this
paper. However, if ${\boldsymbol \Sigma}(\omega) =
\sum_{r\in \mathbb{Z}}{\bf C}(r)e^{ir\omega}$ is a positive definite
matrix with Vector MA$(\infty)$ and Vector AR$(\infty)$
representations (See, \cite{p:wie-mas-58})
 then it  may be possible to extend the above results
to the multivariate setting. 
\end{remark}

\section{Finite order autoregressive approximations}\label{sec:approx}

In many applications, it is often assumed the spectral density is
rational (\cite{p:cad-82,p:ahl-91}, and \cite{p:ge-16}). Obtaining the spectral factorization (such as that given in (\ref{eq:phipsidef})) of a rational spectral density is straightforward, and is one of the reasons that rational spectral densities are widely used.
However, a rational spectral density is usually only an approximation of the underlying spectral
density. In this section, we obtain a bound for the approximation when
the rational spectral density corresponds to a finite order autoregressive
process. The expressions in (\ref{eq:HH}) and (\ref{eq:HH2}) easily
lend themselves to obtaining a rational approximation. Further one can use 
Baxter's inequality to obtain a bound for the approximation. 


We now use the expressions in (\ref{eq:HH})
to obtain an approximation of $H(\omega)$ in terms of the best fitting AR$(p)$ coefficients. 
In particular, using that $\psi(\omega)^* = [\phi(\omega)^*]^{-1}$, we replace the infinite order AR coefficients in
\begin{equation*}
H(\omega) =
\frac{\sum_{\ell=0}^{\infty}
g_{\ell} \big( e^{i\ell\omega} + [\phi(\omega)^*]^{-1} \phi_{\ell}(\omega)^*\big)}{f_{}(\omega)}
\end{equation*}
with the best fitting AR$(p)$ coefficients. 
More precisely, suppose that $(\phi_{p,1}, ..., \phi_{p,p})^\prime$ are 
the best fitting AR$(p)$ coefficients in the sense that it
minimizes the mean squared prediction error 
\begin{equation}
\label{eq:32fit}
(\phi_{p,1}, ..., \phi_{p,p})^\prime =
 \arg\min_{\bf a} \Ex \big| X_{0} - \sum_{j=1}^{p} a_{j}X_{-j}\big|^2 =
 \arg\min_{\bf a} \frac{1}{2\pi}
\int_{0}^{2\pi}\big|1-\sum_{j=1}^{p}a_{j}e^{ij\omega}\big|^{2}f(\omega)d\omega,
\end{equation} 
where ${\bf a} = (a_1, ..., a_p)^\prime$. The corresponding best fitting AR$(p)$ spectral density is $f_{p}(\omega) = \sigma_p^{2}|\phi^{(p)}(\omega)|^{-2}$ where
$\sigma_p^{2} = \Ex |X_0 -\sum_{j=1}^{p} \phi_{p,j} X_{-j}|^2$ and $\phi^{(p)}(\omega) = 1-\sum_{j=1}^{p}\phi_{p,j}e^{ij\omega}$. 
We note that the zeros of the characteristic polynomial $1-\sum_{j=1}^{p}\phi_{p,j}z^j$ lie outside the unit circle (see \cite{b:bro-dav-06}, Problem 8.3).
Then, we define the approximation of $H(\omega)$ as
\begin{equation}
\label{eq:S2}
H_{p}(\omega) = 
\frac{\sum_{\ell=0}^{\infty}
g_{\ell} \big( e^{i\ell\omega} +    [\phi^{(p)}(\omega)^*]^{-1} \phi_{\ell}^{(p)}(\omega)^*\big)}{f_{p}(\omega)},
\end{equation} where $\phi_{\ell}^{(p)}(\omega) = \sum_{s=1}^{p-\ell} \phi_{p,\ell+s}e^{is\omega}$ for $0\leq \ell<p$ and 0 for $\ell \geq p$.
We observe that the Fourier coefficients of $H_{p}(\omega)$ are the solution of 
$T(f_{p}) {\bf h}_{p} = {\bf g}_+$ where ${\bf h}_{p} = (h_{p,0},h_{p,1},\ldots)^\prime$ with
$h_{p,j}
=(2\pi)^{-1}\int_{0}^{2\pi}H_{p}(\omega)e^{-ij\omega}d\omega$. Thus 
$T(f_{p})$ and $T(f_{p})^{-1}$ are approximations of $T(f)$ and
$T(f)^{-1}$ respectively. By using Lemma \ref{lemma:verification} and 
(\ref{eq:HHPP}) we can show that
\begin{equation}
\label{eq:HHPP2}
H_{p}(\omega) = \sigma_p^{-2}\phi^{(p)}(\omega) 
[ \phi^{(p)}(\omega)^{*} G_{+}(\omega)]_{+}
\end{equation}

From a practical perspective, 
the best fitting AR$(p)$ coefficients
 can be estimated from the data.  The AR coefficients
 $\{\phi_{p,j} : 1\leq j \leq p\}$ in (\ref{eq:HHPP2}) can be replaced by its
 estimate and the result used as an estimator of $H(\omega)$.

Below we obtain a bound for $H(\omega)-H_{p}(\omega)$. 

\begin{theorem}[Approximation theorem] \label{theorem:approx}
Suppose that $\{c(r):r\in\mathbb{Z}\}$ is a symmetric, positive definite
sequence that satisfies Assumption \ref{assum:A}(ii) 
and its Fourier transform $f(\omega) = \sum_{r\in
  \mathbb{Z}}c(r)e^{ir\omega}$ 
satisfies Assumption \ref{assum:A}(i).  We define the (semi) infinite system of equations 
\begin{equation*}
g_{\ell} = \sum_{j=0}^{\infty}h_{j}c(\ell-j) \qquad \ell \geq 0,
\end{equation*}
where $(g_{\ell}; \ell \geq 0)^\prime \in \ell_2^+$.  Let $H(\omega)$ and $H_{p}(\omega)$ be defined as in (\ref{eq:Homega}) and (\ref{eq:S2}). Then
\begin{equation*}
\big|H(\omega) - H_{p}(\omega)\big| \leq
C\left[p^{-K+1}\sup_{s}|g_{s}| +  p^{-K}|G_+(\omega)|\right],
\end{equation*} 
where $G_+(\omega)=\sum_{\ell=0}^{\infty}g_{\ell}e^{i\ell\omega}$.
\end{theorem}
\noindent PROOF. See Appendix \ref{sec:proof}. \hfill $\Box$

\begin{remark}[Alternative approximation methods]
There are other ways to obtain the spectral
factorization for non-rational spectral density. 
For example, using the Fourier coefficients of the log spectral density $\log
f(\omega)$ (usually called the cepstral coefficients), 
\cite{p:pou-84} proposed a recursive algorithm for obtaining the AR
and MA coefficients (see, also, \cite{p:bau-55}, \cite{p:wil-72} and  \cite{b:mce-19}, Chapter
7.7, Fact 7.7.6). As this is an infinite recursion based on
an infinite number of cepstral coefficients, typically the number
of non-zero cepstral coefficients is truncated to a finite number in
order to terminate the recursion.  The truncation will  lead to an
approximation error, which we do not investigate here. 
\end{remark}



\subsection*{Data Availability Statement}

Data sharing is not applicable to this article as no new data was created or analyzed in this study.

\subsection*{Acknowledgements}

SSR and JY gratefully acknowledge the partial support of the National
Science Foundation (grant DMS-1812054). JY's research was also supported from the Ministry of Science and Technology, Taiwan (grant 110-2118-M-001-014-MY3).
The authors are extremely gratefully to the comments and corrections
made by two anonymous referees. Their insights substantially improved all aspects of the paper.

\appendix

\section{Proofs} \label{sec:proof}

The purpose of this appendix is to give the technical details behind the results
stated in the main section. 

\vspace{1em}

\noindent {\bf PROOF of Theorem \ref{theorem:wiener}} \hspace{0.3em}
To prove that ${\bf h}_+ = (h_j; j\geq 0)^\prime \in \ell_{2}^{+}$, we note that
since $\Ex[Y^{2}]<\infty$, then $P_{\mathcal{H}_0}(Y) = \sum_{j=0}^{\infty}h_{j}X_{-j}$ is a well-defined random variable in $\mathcal{H}_0$ with
\begin{equation*}
\var[P_{\mathcal{H}_0}(Y)] = \langle {\bf h}_{+},T(f) {\bf h}_{+} \rangle < \infty.
\end{equation*} 
Furthermore, we note that
\begin{equation*}
\langle {\bf h}_{+},T(f) {\bf h}_{+} \rangle =
 \frac{1}{2\pi} \int_{0}^{2\pi} \big| \sum_{j=0}^{\infty} h_j e^{ij\omega}\big|^2 f(\omega) d\omega \geq 
\inf_{\omega} f(\omega) \cdot
\frac{1}{2\pi} \int_{0}^{2\pi} \big| \sum_{j=0}^{\infty} h_j e^{ij\omega}\big|^2 d\omega.
\end{equation*}  
Since $\inf_\omega f(\omega) >0$, we have $\sum_{j=0}^{\infty} h_j e^{ij\omega} \in L_2([0,2\pi))$ and thus
${\bf h}_{+} \in \ell_2^+$.

To prove that
${\bf c}_{+} =(c_{YX}(\ell); \ell \geq 0)^{\prime} \in \ell_{2}^{+}$, we
recall that (\ref{eq:PY}) leads to the matrix equation
${\bf c}_{+} = T(f){\bf h}_+$.
Let $\|A\|_{sp} = \sup_{v \in \ell_2^+, \|v\|_2 = 1} \|Av\|_2$ be the spectral norm. Then, since $\sup_\omega f(\omega) <\infty$,
$\|T(f)\|_{sp} \leq \sup_\omega f(\omega) <\infty$, we have that 
\begin{equation*}
\|{\bf c}_{+}\|_2 = \|T(f){\bf h}_+\|_2 \leq \|T(f)\|_{sp} \|{\bf h}_+\|_2 <\infty.
\end{equation*} 
Thus, ${\bf c}_{+} \in \ell_2^+$.

From (\ref{eq:HFF}), we have 
$H(\omega) = F({\bf c}_{\pm})(\omega)/f_{}(\omega)$. Our goal is to express $F({\bf c}_{\pm})(\omega)$ in terms of the infinite order 
AR and MA coefficients of $\{X_{t}\}$. To do this we observe 
\begin{equation}
\label{eq:Fcpn}
F({\bf c}_{\pm})(\omega)= \sum_{\ell =
  0}^{\infty}c_{YX}(\ell)e^{i\ell\omega} 
+\sum_{\ell=1}^{\infty}\cov\big( Y,P_{\mathcal{H}_0}(X_{\ell}) \big)e^{-i\ell\omega}.
\end{equation}
The second term on the right hand side of (\ref{eq:Fcpn}) looks quite
unwieldy. However, we show below that it can be expressed in terms of the infinite order
AR coefficients associated with $f$. It is well-known that the $\ell$-step ahead forecast $P_{\mathcal{H}_0}(X_{\ell})$
($\ell >0$) has the representation 
$P_{\mathcal{H}_0}(X_{\ell}) = \sum_{j=0}^{\infty}\phi_{j}(\ell)X_{-j}$ with $\ell$-step ahead prediction coefficients
\begin{equation}
\label{eq:phiellj}
\phi_{j}(\ell) = \sum_{s=1}^{\ell}\phi_{j+s}\psi_{\ell-s},
\end{equation}
where $\{\phi_{j}:j\geq 1\}$ and $\{\psi_j: j\geq 0\}$ are the infinite order AR
and MA coefficients defined in (\ref{eq:phipsidef}) (setting
$\psi_0=1$), respectively. 
We now obtain an expression for 
$\cov( Y,P_{\mathcal{H}_0}(X_{\ell}) )$. Using (\ref{eq:phiellj}),
\begin{equation}
\label{eq:PredZ}
 \cov( Y,P_{\mathcal{H}_0}(X_{\ell}))= \cov\big( Y, \sum_{j=0}^{\infty} \sum_{s=1}^{\ell}\phi_{j+s}\psi_{\ell-s} X_{-j}   \big) =
\sum_{s=1}^{\ell} \sum_{j=0}^{\infty}c_{YX}(j) \phi_{j+s}\psi_{\ell-s}.
\end{equation}
For the second identity above, we use Fubini's theorem; noting that
coefficients are absolutely summable since 
\begin{equation*}
\sum_{s=1}^{\ell} |\psi_{\ell-s}| \sum_{j=0}^{\infty}   |c_{YX}(j) \phi_{j+s}|
\leq \left( \sum_{s=1}^{\ell}|\psi_{\ell-s}| \right) 
\left( \sum_{j=0}^{\infty} c(j)^2 \right)^{1/2} \left( \sum_{j=1}^{\infty} \phi_j^2\right)^{1/2}
<\infty.
\end{equation*}
Using (\ref{eq:PredZ}) we have 
\begin{equation*}
\sum_{\ell=1}^{\infty}\cov(Y,P_{\mathcal{H}_0}(X_{\ell}))
  e^{-i\ell\omega} = \sum_{\ell=1}^{\infty} \left( \sum_{s=1}^{\ell} 
\psi_{\ell-s} 
\sum_{j=0}^{\infty}c_{YX}(j) \phi_{j+s}\right) e^{-i\ell\omega}.
\end{equation*} 
The Fourier coefficients of the right hand side of above has a convolution form,
thus, we use the convolution theorem and rewrite
\begin{eqnarray*}
\sum_{\ell=1}^{\infty}\cov(Y,P_{\mathcal{H}_0}(X_{\ell}))
  e^{-i\ell\omega} &=& \sum_{\ell=1}^{\infty} \left( \sum_{s=1}^{\ell} 
\psi_{\ell-s} 
\sum_{j=0}^{\infty}c_{YX}(j) \phi_{j+s}\right) e^{-i\ell\omega} \\
&=&
\left( \sum_{\ell=0}^{\infty} \psi_\ell e^{-i\ell \omega}\right)
\left( \sum_{s=1}^{\infty} \sum_{j=0}^{\infty} c_{YX}(j) \phi_{j+s} e^{-is\omega}\right) \\
&=& \psi(\omega)^* \sum_{j=0}^{\infty} c_{YX}(j)
\sum_{s=1}^{\infty}  \phi_{j+s} e^{-is\omega}  \\
&=& \psi(\omega)^* \sum_{j=0}^{\infty} c_{YX}(j) \phi_j(\omega)^*,
\end{eqnarray*} 
where $\phi_{j}(\omega) = \sum_{s=1}^{\infty}\phi_{j+s}e^{is\omega}$ for $j \geq 0$.
Substituting the above into (\ref{eq:Fcpn}) gives 
\begin{eqnarray}
F({\bf c}_{\pm})(\omega) &=& \sum_{\ell =
  0}^{\infty}c_{YX}(\ell)e^{i\ell\omega}  +
\psi_{}(\omega)^{*} \sum_{j=0}^{\infty} \phi_{j}(\omega)^{*} c_{YX}(j) \nonumber\\
&=&  \sum_{\ell =
  0}^{\infty}c_{YX}(\ell) \big( e^{i\ell\omega} +
    \psi_{}(\omega)^*\phi_{\ell}(\omega)^*\big).
 \label{eq:atilde}
\end{eqnarray}
Since $\psi_{}(\omega)^*$ is bounded, it is easily seen that $F({\bf c}_{\pm})(\omega)\in
L_{2}([0,2\pi))$. Finally, substituting the above into  $H(\omega)=F({\bf c}_{\pm})(\omega)/f_{}(\omega)$ proves (\ref{eq:Homega}).

To prove the alternative expression in (\ref{eq:alternativeform}), we
rearrange the expression $e^{i\ell \omega} + \psi(\omega)^*
\phi_\ell(\omega)^*$ which appears in $F({\bf
  c}_{\pm})(\omega)$. Using the definition $\phi(\omega)^* =
1-\sum_{j=1}^{\infty} \phi_je^{-ij\omega}$, we have 
\begin{eqnarray}
e^{i\ell \omega} + \psi(\omega)^* \phi_\ell(\omega)^* &=&
e^{i\ell \omega} ( 1+ \psi(\omega)^* \sum_{s=1}^{\infty} \phi_{s+\ell}e^{-i(s+\ell)\omega}) \nonumber \\
&=&
e^{i\ell \omega} \left( 1+ \psi(\omega)^* \left[-\phi(\omega)^* +1 -
    \sum_{j=1}^{\ell} \phi_j e^{-ij\omega}\right] \right) 
\nonumber \\
&=& e^{i\ell \omega}\psi(\omega)^* \bigg( 1-\sum_{j=1}^{\ell} \phi_j
    e^{-ij\omega} \bigg). 
\qquad \text{(using $\phi(\omega)^* \psi(\omega)^* = 1$)}
\label{eq:alt-exp}
\end{eqnarray} 
Therefore, substituting (\ref{eq:alt-exp}) into
(\ref{eq:atilde}) and using that 
$H(\omega)=F({\bf c}_{\pm})(\omega)/f_{}(\omega)$ gives
\begin{eqnarray*}
H(\omega)= \frac{F({\bf c}_{\pm})(\omega)}{f(\omega)}
&=& \frac{\sum_{\ell =0}^{\infty}c_{YX}(\ell) e^{i\ell \omega}\psi(\omega)^* \bigg( 1-\sum_{j=0}^{\ell} \phi_j e^{-ij\omega} \bigg)
}{f(\omega)} \\
&=& \frac{\sum_{\ell =0}^{\infty}c_{YX}(\ell) \bigg( e^{i\ell \omega} - \sum_{j=1}^{\ell} \phi_j e^{i(\ell-j)\omega} \bigg)
}{\sigma^2 \psi(\omega)} \quad \text{(using $f(\omega) = \sigma^2 \psi(\omega) \psi(\omega)^*$)} \\
&=& \sigma^{-2} \phi(\omega) \sum_{\ell =0}^{\infty}c_{YX}(\ell) \bigg( e^{i\ell \omega} - \sum_{j=1}^{\ell} \phi_j e^{i(\ell-j)\omega} \bigg).
\end{eqnarray*} This shows (\ref{eq:alternativeform}) and thus proves the Theorem.
\hfill $\Box$

\vspace{1em}

\noindent \textbf{Proof of equation (\ref{eq:PPY}) in  Remark \ref{remark:wiener}}  \hspace{0.3em}
For fixed $j \geq 0$, the coefficient of $X_{-j}$ in 
$\sum_{\ell=1}^{\infty}P_{\mathcal{H}_0}(X_{\ell})e^{i\ell\omega}$ is $\sum_{\ell=1}^{\infty} \phi_j(\ell) e^{i\ell\omega}$. Using (\ref{eq:phiellj}), we get
\begin{eqnarray*}
\sum_{\ell=1}^{\infty} \phi_j(\ell) e^{i\ell\omega} &=& \sum_{\ell=1}^{\infty} 
\sum_{s=1}^{\ell}\phi_{j+s}\psi_{\ell -s} e^{i\ell\omega} \\
&=&
\sum_{s=1}^{\infty}\phi_{j+s} e^{is\omega}
\sum_{\ell=s}^{\infty}\psi_{\ell -s}
e^{i(\ell-s)\omega} \\
&=&
    \psi(\omega) \sum_{s=1}^{\infty}\phi_{j+s}e^{is\omega}= 
 \psi(\omega)\phi_{j}(\omega).
\end{eqnarray*} The second identity above is also due to the fact that $\sum_{s=1}^{\infty} \sum_{\ell=s}^{\infty} |\phi_{j+s}\psi_{\ell -s}| < \infty$ under Assumption \ref{assum:A}(ii) for $K=1$.
We now show that under Assumption \ref{assum:A}(ii) for $K=1$, 
$\psi(\omega)\sum_{j=0}^{\infty}X_{-j}\phi_{j}(\omega)$ converges in
$\mathcal{H}_0$.  To show this, we
define the partial sum
\begin{equation*}
S_{n} = \psi(\omega)\sum_{j=0}^{n}X_{-j}\phi_{j}(\omega) \in \mathcal{H}_0.
\end{equation*}
Then, for any $n<m$
\begin{equation*}
\Ex |S_{m}-S_{n}|^{2} = 
\var \left( \sum_{j=n}^{m} \psi(\omega) \phi_j(\omega) X_{-j} \right)
= |\psi(\omega)|^2 ({\bf \phi}_{n}^{m}(\omega))^\prime T_{m-n}(f) ({\bf \phi}_{n}^{m}(\omega)),
\end{equation*} 
where $({\bf \phi}_{n}^{m}(\omega)) = (\phi_n(\omega), ...,
\phi_{m}(\omega))^\prime$ and 
$T_{m-n}(f) = (c(t-\tau); 0 \leq t, \tau \leq m-n)$. Therefore,
\begin{eqnarray*}
\Ex |S_{m}-S_{n}|^{2}
&=& |\psi(\omega)|^2 ({\bf \phi}_{n}^{m}(\omega))^\prime T_{m-n}(f) ({\bf \phi}_{n}^{m}(\omega)) \\
&\leq& |\psi(\omega)|^2 \|T_{m-n}(f)\|_{spec} \|{\bf \phi}_{n}^{m}(\omega)\|_2^2  \\
&\leq& |\psi(\omega)|^2 ( \sup_\omega f(\omega) )\|{\bf \phi}_{n}^{m}(\omega)\|_2^2.
\end{eqnarray*}
If Assumption \ref{assum:A}(ii) is satisfied for $K=1$, then it is easy to show
$\sum_{j=0}^{\infty} |\phi_j(\omega)|^2 <\infty$.
Therefore, by Cauchy's criterion, 
$\|{\bf \phi}_{n}^{m}(\omega)\|_2 \rightarrow 0$ as $n,m \rightarrow \infty$, which implies $\Ex |S_{m}-S_{n}|^{2} \rightarrow 0$ as $n,m \rightarrow \infty$.
Again applying Cauchy's criterion (on the Hilbert space $\mathcal{H}_0$), we conclude that $\psi(\omega)\sum_{j=0}^{\infty}X_{-j}\phi_{j}(\omega)$ converges in
$\mathcal{H}_0$. 
This shows $\sum_{\ell=1}^{\infty}P_{\mathcal{H}_0}(X_{\ell}) e^{i\ell\omega}$ is well-defined in $\mathcal{H}_0$ and
satisfies (\ref{eq:PPY}).
\hfill $\Box$

\vspace{1em}

\noindent \textbf{PROOF of equation (\ref{eq:C01C11})}
\hspace{0.3em}
Representing 
\begin{equation}
\label{eq:PH0}
P_{\mathcal{H}_0}(X_{-\ell}) = \sum_{j=0}^{\infty} A_{\ell,j} X_{-j} \qquad \ell <0
\end{equation}
we will show that matrix $A = (A_{\ell,j}; \ell <0, j\geq 0) = C_{01}C_{11}^{-1}$. We
first evaluate the covariance $\cov(X_{-t}, P_{\mathcal{H}_0}(X_{-\ell}))$
 (for $t \geq 0$) using (\ref{eq:PH0}). The 
left hand side of (\ref{eq:PH0}) is
\begin{equation*}
\cov(X_{-t}, P_{\mathcal{H}_0}(X_{-\ell}))
 = \cov(P_{\mathcal{H}_0}(X_{-t}), X_{-\ell})
 = c(\ell-t) = [C_{01}]_{\ell,t} \quad \ell<0, t\geq 0.
\end{equation*} 
Whereas the right hand side of (\ref{eq:PH0}) is
\begin{equation*}
\cov (X_{-t}, \sum_{j=0}^{\infty} A_{\ell,j} X_{-j})
= \sum_{j=0}^{\infty} A_{\ell,j} c(j-t) = \sum_{j=0}^{\infty} A_{\ell,j} [C_{11}]_{j,t} = [A C_{11}]_{\ell,t} \quad \ell<0, t\geq 0.
\end{equation*} 
Comparing coefficients gives $C_{01} = AC_{11}$, i.e., $A = C_{01}C_{11}^{-1}$.
\hfill $\Box$

\vspace{1em}

\noindent \textbf{PROOF of Theorem \ref{corollary:wiener}}
\hspace{0.3em}
We first prove that ${\bf h}_+ = (h_0, h_1, ...)^\prime \in \ell_{2}^{+}$. 
Under Assumption \ref{assum:A}(i), $T(f)$ is invertible on $\ell_2^+$ (\cite{p:wid-60}, Theorem III).
Using that $\|T(f)^{-1}\|_{sp} \leq [\inf_{\omega}f(\omega)]^{-1}$ we
have 
\begin{equation*}
\|{\bf h}_+\|_{2}\leq \|T(f)^{-1}\|_{sp}\|{\bf g}_{+}\|_{2}\leq
  [\inf_{\omega}f(\omega)]^{-1}\|{\bf g}_{+}\|_{2},
\end{equation*} where ${\bf g}_{+} = (g_0, g_1, ...)^\prime$.
Since ${\bf g}_{+} \in \ell_{2}^{+}$, from the above inequality, we get
${\bf h}_{+} \in \ell_{2}^{+}$ and its Fourier transform $H(\omega)$ is well-defined.
Thus, using the construction described as in (\ref{eq:Ydef}),
there exists a second order stationary time series
$\{X_{t}\}$ and random variable $Y \in \overline{\spa} (X_t: t\leq 0)$,
whose normal equations satisfy
\begin{equation*}
\cov (Y, X_{-\ell}) = g_{\ell}  =
  \sum_{j=0}^{\infty}h_{j}c(\ell-j) \qquad \ell \geq 0.
\end{equation*}
This allows us to use  Theorem \ref{theorem:wiener} to prove the
result. 
\hfill $\Box$

\vspace{1em}

\noindent {\bf PROOF of Lemma \ref{lemma:verification}} \hspace{0.3em}
Using that $\phi(\omega)^*= 1-\sum_{j=1}^{\infty} \phi_j e^{-ij\omega}$, the right hand side of (\ref{eq:link42}) is
\begin{eqnarray*}
\sum_{\ell=0}^{\infty}
g_\ell \left( e^{i\ell \omega} - \sum_{s=1}^{\ell} \phi_s e^{i(\ell-s)\omega} \right) &=&
\sum_{\ell=0}^{\infty} g_\ell e^{i\ell\omega} \left( 1 - \sum_{s=1}^{\ell} \phi_s e^{-is\omega}\right) \\
&=&
\sum_{\ell=0}^{\infty} g_\ell e^{i\ell\omega} \left( \phi(\omega)^* + \sum_{s=\ell+1}^{\infty} \phi_s e^{-is\omega} \right) \\
&=& G_{+}(\omega) \phi(\omega)^* + 
\sum_{\ell=0}^{\infty} g_\ell \sum_{s=\ell+1}^{\infty} \phi_s e^{i(\ell-s)\omega}.
\end{eqnarray*} It is straightforward that the second term on the right hand side above is anti-causal. Therefore, take the causal part of above gives
\begin{equation*}
\left[ \sum_{\ell=0}^{\infty}
g_\ell \left( e^{i\ell \omega} - \sum_{s=1}^{\ell} \phi_s e^{i(\ell-s)\omega} \right)\right]_{+}
= [ G_{+}(\omega) \phi(\omega)^*]_+.
\end{equation*}
 This proves the lemma. 
\hfill $\Box$

\vspace{1em}

\noindent {\bf Proof of Corollary \ref{theorem:inverse}} \hspace{0.3em}
Let
$\delta_{\ell,k}$ denotes the indicator variable where
$\delta_{\ell,k}=1$ if $\ell=k$ and zero otherwise.
Since $T(f)^{-1} = (d_{j,k}; j,k \geq 0)$ is the inverse of $T(f) = (c(j-k); j,k \geq 0)$,
 $\{d_{j,k}\}$ and $\{c(r)\}$ satisfy the normal equations
\begin{equation}
\label{eq:deltajell}
\delta_{\ell,k} = \sum_{j=0}^{\infty} d_{j,k} c(\ell-j) 
\qquad \ell,k \geq 0.
\end{equation}
Thus for each $k \geq 0$, we obtain a system of Wiener-Hopf equations. 
To derive $d_{j,k}$ we apply  Theorem
\ref{corollary:wiener} to (\ref{eq:deltajell}).  For each (fixed) $k\geq 0$ we obtain
\begin{equation}
\label{eq:Dkomega1}
D_{k}(\omega) =
\frac{1}{f(\omega)}\sum_{\ell =0}^{\infty}
\delta_{\ell,k} \big( e^{i\ell\omega} +
    \psi_{}(\omega)^*\phi_{\ell}(\omega)^*\big)
= \frac{e^{ik\omega} +
    \psi_{}(\omega)^*\phi_{k}(\omega)^*}{f(\omega)},
\end{equation}
where $D_{k}(\omega) = \sum_{j=0}^{\infty}d_{j,k}e^{ij\omega}$. Using
the identity (\ref{eq:alt-exp})
we can replace the above with
\begin{equation} \label{eq:Dkomega2}
D_{k}(\omega)
= \frac{\psi(\omega)^* (e^{ik\omega} + \sum_{s=1}^{k} \phi_{s}e^{i(k-s)\omega})}{f(\omega)}
= \sigma^{-2} \phi(\omega) \bigg(e^{ik\omega} + \sum_{s=1}^{k} \phi_{s}e^{i(k-s)\omega}\bigg).
\end{equation}
Taking an inverse Fourier transform in (\ref{eq:Dkomega1}) and (\ref{eq:Dkomega2}) yields the entries
\begin{eqnarray*}
d_{k,j}  &=& \frac{1}{2\pi}\int_{0}^{2\pi}\left(\frac{e^{ik\omega} +
    \psi_{}(\omega)^*\phi_{k}(\omega)^*}{f(\omega)}\right)e^{-ij\omega}d\omega \\
&=& \sigma^{-2} \frac{1}{2\pi}\int_{0}^{2\pi}
\phi(\omega) \bigg(e^{ik\omega} + \sum_{s=1}^{k} \phi_{s}e^{i(k-s)\omega}\bigg) e^{-ij\omega} d\omega \qquad j,k \leq 0.
\end{eqnarray*}
Thus proving the Corollary. 

As an aside it  is interesting to construct the random variable
$Y_{-k}$ which yields the Wiener-Hopf equation (\ref{eq:deltajell}).
Let $\{X_{t}\}$ be a second order stationary process with
$\{c(r)\}$ as its autocovariance. We define a sequence of random variables $\{\varepsilon_{-k}: k \geq 0\}$
where for $k \geq 0$
\begin{equation*}
\varepsilon_{-k} = X_{-k} -P_{(-k)^c}(X_{-k})
\end{equation*}
and $P_{(-k)^c}$ denotes the orthogonal projection 
onto the closed subspace $\overline{\textrm{sp}}(X_{t};t\leq 0\textrm{ and
}r\neq -k)$. We standardize $\varepsilon_{-k}$, where
$Y_{-k} =
\varepsilon_{-k}/\sqrt{\var(\varepsilon_{-k})}$, noting that
 $\var(\varepsilon_{-k}) = \cov(\varepsilon_{-k},X_{-k})$. Thus by
 definition $\cov(Y_{-k},X_{\ell}) = \delta_{\ell,k}$ and 
$Y_{-k} = \sum_{j=0}^{\infty} d_{k,j}X_{-j}$.
\hfill $\Box$

\vspace{1em}

\noindent {\bf PROOF of Theorem \ref{theorem:approx}} \hspace{0.3em}
We note that under Assumption \ref{assum:A}(ii), $\sum_{j=1}^{\infty}|j^{K}\phi_{j}|<\infty$.

To prove the result, we use Baxter's inequality, that is for the best fitting
AR$(p)$ coefficients (see equation (\ref{eq:32fit})), we have 
\begin{equation}
\label{eq:baxter}
\sum_{j=1}^{p}|\phi_{p,j}-\phi_{j}| \leq C_{f} \sum_{j=p+1}^{\infty}|\phi_{j}|
\end{equation}
where $C_{f}$ is a constant that soley depends on $f(\omega) = \sigma^{2}|\phi(\omega)|^{-2}$.

\vspace{1em}

Returning to the proof, the difference $H(\omega) - H_{p}(\omega)$ can
be decomposed as
\begin{eqnarray*}
H(\omega) - H_{p}(\omega) &=&
  \sum_{\ell=0}^{\infty}g_{\ell}e^{i\ell\omega}\left(\frac{1}{f(\omega)}
  - \frac{1}{f_{p}(\omega)}\right) 
 + \sum_{\ell = 0}^{\infty}g_{\ell} 
\left(\frac{[\phi(\omega)^*]^{-1}\phi_{\ell}(\omega)^*}{f(\omega)} -
   \frac{[\phi^{(p)}(\omega)^*]^{-1}\phi_{\ell}^{(p)}(\omega)^*}{f_{p}(\omega)}   \right) \\
&=&
  \sum_{\ell=0}^{\infty}g_{\ell}e^{i\ell\omega}\left(\frac{1}{f(\omega)}
  - \frac{1}{f_{p}(\omega)}\right) 
 + \sum_{\ell = 0}^{\infty}g_{\ell} 
\big[ \sigma^{-2}\phi(\omega) \phi_{\ell}(\omega)^*-
   \sigma_p^{-2}\phi^{(p)}(\omega) \phi_{\ell}^{(p)}(\omega)^* \big] \\
&=& A(\omega) + B(\omega)+ C(\omega)
\end{eqnarray*} 
where
\begin{eqnarray*}
A(\omega) &=& \left(\frac{1}{f(\omega)} - \frac{1}{f_{p}(\omega)}
      \right)\sum_{\ell=0}^{\infty}g_{\ell}e^{i\ell\omega} \\
B(\omega) &=& \left[\sigma^{-2}\phi(\omega)-\sigma_{p}^{-2}\phi^{(p)}(\omega)
      \right]\sum_{\ell=0}^{\infty}g_{\ell}\phi_{\ell}(\omega)^* \\
C(\omega) &=& \sigma_p^{-2}
      \phi^{(p)}(\omega)\sum_{\ell=0}^{\infty}g_{\ell}\left[\phi_{\ell}(\omega)^*
      - \phi_{\ell}^{(p)}(\omega)^*\right].
\end{eqnarray*}
To bound each of the terms above, we derive some auxillary bounds
which we use later to bound $A(\omega)$ to $C(\omega)$.
We use $C$ to denote a generic constant which may change from line to
line (and depends on $C_{f}$).

First, we bound the difference $|\phi(\omega) - \phi^{(p)}(\omega)|$. Using (\ref{eq:baxter}), we have
\begin{eqnarray}
|\phi(\omega)- \phi^{(p)}(\omega)| &=& |\sum_{j=1}^{p} (\phi_j - \phi_{p,j})e^{ij\omega}
+ \sum_{j=p+1}^{\infty} \phi_j e^{ij\omega} | \nonumber \\
 &\leq& \sum_{j=1}^{p} |\phi_j -\phi_{p,j}| + \sum_{j=p+1}^{\infty} |\phi_j| 
\leq (C_{f} + 1) \sum_{j=p+1}^{\infty} |\phi_j|.
\label{eq:Phibound}
\end{eqnarray}  
Furthermore, using that $|\phi(\omega)| \leq 1+\sum_{j=1}^{\infty} |\phi_j| < \infty$, we have
\begin{equation*}
|\phi^{(p)}(\omega)| \leq |\phi^{(p)}(\omega) - \phi(\omega)| + |\phi(\omega)|
\leq (C_{f}+1) \sum_{j=p+1}^{\infty} |\phi_j| + 1+\sum_{j=1}^{\infty} |\phi_j| < \infty.
\end{equation*} 
Using above two bounds, we obtain a bound for 
$|\phi(\omega)|^2 - |\phi^{(p)}(\omega)|^2 $. 
By the triangular inequalities $||A^2| -|B^2|| \leq |A^2-B^2|$ and
$|A+B|\leq |A|+|B|$ together with (\ref{eq:Phibound}) gives 
\begin{eqnarray}
| |\phi(\omega)|^2 - |\phi^{(p)}(\omega)|^2 | \leq
| \phi(\omega)^2 - \phi^{(p)}(\omega)^2 | &\leq&
|\phi(\omega) - \phi^{(p)}(\omega)| (|\phi(\omega)| + |\phi^{(p)}(\omega)|) \nonumber \\
&\leq& C \sum_{j=p+1}^{\infty} |\phi_j|.
\label{eq:Phibound1}
\end{eqnarray}
Next, we bound $|\sigma^{-2} - \sigma_p^{-2}|$. 
We recall that $\sigma^{2} = (2\pi)^{-1}\int
  |\phi(\omega)|^{2}f(\omega)d\omega$ and 
$\sigma^{2}_{p} =
  \Ex|X_{t}-\sum_{j=1}^{p}\phi_{p,j}X_{t-j}|^{2} = (2\pi)^{-1}\int
  |\phi_{p}(\omega)|^{2}f(\omega)d\omega$. Using these expression 
we have
\begin{equation*}
|\sigma^2 - \sigma_p^2| \leq \frac{1}{2\pi} \int_{0}^{2\pi}
| |\phi(\omega)|^2 - |\phi^{(p)}(\omega)|^2 | f(\omega) d\omega. 
\end{equation*} 
Combining the above with (\ref{eq:Phibound1}) 
and Assumption \ref{assum:A}(i)
gives
\begin{eqnarray*}
|\sigma^2 - \sigma_p^2| &\leq& \frac{1}{2\pi} \int_{0}^{2\pi}
| |\phi(\omega)|^2 - |\phi^{(p)}(\omega)|^2 | f(\omega) d\omega \nonumber \\
&\leq& \sup_\omega f(\omega) \cdot \frac{1}{2\pi} \int_{0}^{2\pi}
| |\phi(\omega)|^2 - |\phi^{(p)}(\omega)|^2 | d\omega 
\leq  C \sum_{j=p+1}^{\infty} |\phi_j|. 
\end{eqnarray*} 
As an immediate consequence of above, we have $\sigma_p^2 > \sigma^2 - |\sigma^2 - \sigma_p^2| \geq \sigma^2/2 > 0$ for large $p$. Since $\sigma_p^2 >0$ for all $p$, we have
\begin{equation} \label{eq:sigmabound1}
\inf_p \sigma_p^2 > 0.
\end{equation} Therefore, we obtain the bound
\begin{equation}
|\sigma^{-2} - \sigma_p^{-2}| \leq \sigma^{-2}\sigma_p^{-2} |\sigma^2 - \sigma_p^2|
\leq C \sum_{j=p+1}^{\infty} |\phi_j|.
\label{eq:sigmabound2}
\end{equation}

\noindent Lastly, by Assumption \ref{assum:A}(ii),
\begin{equation} \label{eq:Phibound2}
\sum_{j=p+1}^{\infty} |j|^{\alpha} |\phi_j| \leq p^{-K+\alpha}
\sum_{j=p+1}^{\infty} |j^K \phi_j| = O(p^{-K+\alpha}) \quad\textrm{
  for }  0\leq \alpha \leq K.
\end{equation}

\vspace{1em}

Now we are ready to bound each term in $H(\omega) - H_p(\omega)$. First, we bound $A(\omega)$. Note that
\begin{eqnarray*}
\left| (f(\omega))^{-1} - (f_{p}(\omega))^{-1} \right|
&=& |\sigma^{-2} |\phi(\omega)|^2 - \sigma_p^2 |\phi^{(p)}(\omega)|^2 | \\
&\leq& \sigma^{-2} ||\phi(\omega)|^2 - |\phi^{(p)}(\omega)|^2|
+|\sigma^{-2} - \sigma_p^{-2}| |\phi^{(p)}(\omega)|^2.
\end{eqnarray*} Using (\ref{eq:Phibound1}) and (\ref{eq:sigmabound2}), we get
\begin{equation*}
\left|(f(\omega))^{-1} - (f_{p}(\omega))^{-1}
      \right|  \leq
\sigma^{-2} ||\phi(\omega)|^2 - |\phi^{(p)}(\omega)|^2|
+|\sigma^{-2} - \sigma_p^{-2}| |\phi^{(p)}(\omega)|^2  \leq C \sum_{j=p+1}^{\infty} |\phi_j|.
\end{equation*}
Therefore, substituting (\ref{eq:Phibound2}) (for $\alpha=0$) into $A(\cdot)$ gives 
\begin{equation*}
|A(\omega)| \leq \left|(f(\omega))^{-1} - (f_{p}(\omega))^{-1}
      \right|  \big|\sum_{\ell=0}^{\infty}g_{\ell}e^{i\ell\omega}\big| \leq
C \left( \sum_{j=p+1}^{\infty}|\phi_{j}| \right) \cdot
      \big|\sum_{\ell=0}^{\infty}g_{\ell}e^{i\ell\omega}\big|
= O\left( p^{-K} \big| G_+(\omega)\big| \right),
\end{equation*} where $G_{+}(\omega) = \sum_{\ell=0}^{\infty} g_{\ell}e^{i\ell \omega}$.

To bound $B(\omega)$ we note that from 
(\ref{eq:Phibound}) and (\ref{eq:sigmabound2})
\begin{equation*}
|\sigma^{-2} \phi(\omega) - \sigma_p^{-2} \phi^{(p)}(\omega)|
\leq \sigma^{-2} |\phi(\omega) - \phi^{(p)}(\omega)|
+ |\sigma^{-2} - \sigma_p^{-2}| |\phi^{(p)}(\omega)| 
\leq C \sum_{j=p+1}^{\infty} |\phi_j|.
\end{equation*}
Therefore, we have 
\begin{equation*}
|B(\omega)| \leq |\sigma^{-2} \phi(\omega) - \sigma_p^{-2} \phi^{(p)}(\omega)|
 \sum_{\ell=0}^{\infty}|g_{\ell}\phi_{\ell}(\omega)|
  \leq C \left( \sum_{j=p+1}^{\infty} |\phi_j| \right) \sum_{\ell=0}^{\infty}|g_{\ell}\phi_{\ell}(\omega)|.
\end{equation*}
The second summand on the right hand side of the above is bounded with
\begin{equation*}
 \sum_{\ell=0}^{\infty}|g_{\ell}\phi_{\ell}(\omega)|
  \leq
\sum_{\ell=0}^{\infty}|g_{\ell}|\sum_{s=1}^{\infty}|\phi_{\ell+s}|
=\sum_{u=1}^{\infty} \sum_{\ell=0}^{u-1}|g_{\ell}|  |\phi_{u}|
\leq \sup_{s}|g_s| \cdot \sum_{u=1}^{\infty} |u\phi_u|.
\end{equation*}
Thus by using the above two bounds and (\ref{eq:Phibound2}) (for $\alpha=0$), we have 
\begin{equation*}
|B(\omega)| \leq
                             C \left( \sum_{j=p+1}^{\infty}|\phi_{j}| \right)
\cdot
\sup_{s}|g_{s}| \cdot \sum_{u=1}^{\infty}|u\phi_{u}|
                             =O\left(\sup_{s}|g_{s}| \cdot p^{-K}\right).
\end{equation*}
Finally, we obtain a bound for $C(\omega)$. Since
$\phi_{\ell}^{(p)}(\omega) =0$ for $\ell \geq p$ we split 
$C(\omega)$ into two parts: $C(\omega) = C_1(\omega) + C_2(\omega)$ where
\begin{eqnarray*}
C_1(\omega) &=& \sigma_p^{-2} \phi^{(p)}(\omega) \sum_{\ell=0}^{p-1} g_{\ell} [\phi_{\ell}(\omega)^* - \phi_{\ell}^{(p)}(\omega)^*] \\ 
C_2(\omega) &=& \sigma_p^{-2} \phi^{(p)}(\omega) \sum_{\ell=p}^{\infty} g_{\ell} \phi_{\ell}(\omega)^*.
\end{eqnarray*} To bound $C_1(\omega)$, we note by (\ref{eq:baxter}) and (\ref{eq:Phibound2}) (for $\alpha=0$)
\begin{eqnarray*}
|\phi_{\ell}(\omega)^* - \phi_{\ell}^{(p)}(\omega)^*| &=& \sum_{s=1}^{p-\ell} |\phi_{\ell+s} - \phi_{p,\ell+s}| + \sum_{s=p-\ell+1}^{\infty} |\phi_{\ell+s}|\\
&\leq& C\sum_{s=p+1}^{\infty} |\phi_s| = O(p^{-K}) \qquad 0\leq \ell <p.
\end{eqnarray*} 
where the $O(p^{-K})$ bound above is uniform over $0\leq \ell
<p$. Therefore, combining the above with (\ref{eq:sigmabound1}) gives
\begin{eqnarray} 
|C_1(\omega)| &\leq& \sigma_p^{-2} |\phi^{(p)}(\omega)| \sum_{\ell=0}^{p-1} |g_{\ell}| |\phi_{\ell}(\omega)^* - \phi_{\ell}^{(p)}(\omega)^*| \nonumber \\
&\leq& C  \sup_{s} |g_s|  \cdot
\sum_{\ell=0}^{p-1}|\phi_{\ell}(\omega)^* - \phi_{\ell}^{(p)}(\omega)^*|
=
O\left( \sup_{s} |g_s| \cdot p^{-K+1}\right).
\label{eq:C1bound}
\end{eqnarray}
To bound $C_2(\omega)$, we use (\ref{eq:Phibound2}) for $\alpha=1$,
\begin{equation*}
\sum_{\ell=p}^{\infty} |g_{\ell}| |\phi_{\ell}(\omega)| \leq
\sup_s |g_{s}| \sum_{\ell=p}^{\infty} \sum_{s=1}^{\infty}|\phi_{\ell+s}|
\leq \sup_s |g_{s}| \sum_{u=p+1}^{\infty} |u \phi_u|
= O\left( \sup_{s} |g_s| \cdot p^{-K+1}\right).
\end{equation*} Therefore, we have
\begin{equation} \label{eq:C2bound}
|C_2(\omega)| \leq \sigma_p^{-2} |\phi^{(p)}(\omega)| \sum_{\ell=p}^{\infty} |g_{\ell}| |\phi_{\ell}(\omega)| =
O\left( \sup_{s} |g_s| \cdot p^{-K+1}\right).
\end{equation} Combining (\ref{eq:C1bound}) and (\ref{eq:C2bound}) gives
\begin{equation*}
|C(\omega)| \leq |C_1(\omega)| + |C_2(\omega)| = 
O\left( \sup_{s} |g_s| \cdot p^{-K+1}\right).
\end{equation*}
Altogether, this yields the bound 
\begin{equation*}
\big|H(\omega) - H_{p}(\omega)\big| \leq
C\left[ p^{-K+1} \cdot \sup_{s}|g_{s}| + p^{-K} \cdot \big|G_+(\omega)\big|\right].
\end{equation*} 
This proves the result. \hfill $\Box$

\vspace{2mm}

\bibliography{bib_Wiener}
\bibliographystyle{plainnat}

\end{document}